\documentclass[12pt]{article}
\usepackage[T2A]{fontenc}
\usepackage[cp1251]{inputenc}
\usepackage[russian]{babel}
\usepackage{amssymb,amsmath,latexsym,amsthm}
\usepackage{amsbsy,amsfonts,amssymb,mathrsfs,amscd}
\usepackage{eucal}
\usepackage{graphicx}
\usepackage{indentfirst}
\usepackage{verbatim}
\usepackage{enumerate}
\usepackage[pdftex,unicode,colorlinks,linkcolor=blue,
citecolor=red,bookmarksopen,pdfhighlight=/N]{hyperref}

\theoremstyle{plain}

\newtheorem{lemma}{Лемма}
\newtheorem{theorem}{Теорема}

\newtheorem{thA}{Теорема Мальявена\,--\,Рубела\!\!\!}
\newtheorem{thD}{Теорема Kh\!\!\!}

\theoremstyle{definition}

\newtheorem{remark}{Замечание}

\renewcommand{\leq}{\leqslant}
\renewcommand{\geq}{\geqslant}

\newcommand{\rad}{\text{\tiny\rm rad}}
\newcommand{\RR}{\mathbb{R}} 
\newcommand{\CC}{\mathbb{C}} 
\newcommand{\NN}{\mathbb{N}} 

\DeclareMathOperator{\Zero}{\text{\sf Zero}}

\DeclareMathOperator{\type}{\text{\sf type}} 
\DeclareMathOperator{\Hol}{Hol}

\DeclareMathOperator{\rh}{\text{\rm \tiny rh}}
 \DeclareMathOperator{\lh}{\text{\rm \tiny lh}}
\DeclareMathOperator{\dd}{\,{\mathrm d\!}}
\renewcommand{\Re}{{\rm Re \,}}
\renewcommand{\Im}{{\rm Im \,}}

\title{К теореме Мальявена\,--\,Рубела о малости роста целых функций экспоненциального типа с заданными нулями}
\author{А.\,Е.~Салимова, Б.\,Н.~Хабибуллина}
\date{30 апреля 2021 г.}
\begin{document}
\maketitle

Одноточечные множества $\{x\}$ часто записываем без фигурных скобок, т.е. просто как $x$. 
Как обычно, $\mathbb N:=\{1,2, \dots\}$ множество 
{\it натуральных чисел,\/}  $\NN_0:=0\cup \NN$. 
Множество $\RR$ {\it действительных чисел,\/} со стандартными порядковой
($\le$, $\sup/\inf$), алгебраической и топологической структурами в основном рассматривается  как  {\it вещественная ось\/} в  {\it комплексной плоскости\/} $\CC$; $i\RR$ --- {\it мнимая ось},  
$\overline \RR:=-\infty\cup \RR\cup +\infty$ --- {\it расширенная действительная прямая\/} с двумя концами $\pm\infty\notin \RR$, дополненная неравенствами $-\infty\leq x\leq +\infty$ для любого $x\in \overline \RR$ и снабжённая естественной порядковой топологией, а $\overline \NN:=\NN\cup +\infty$, $\overline \NN_0:=\NN_0\cup +\infty$.
{\it Интервалы на\/} $\overline{\RR}$ --- связные подмножества в $\overline \RR$, такие, как  {\it отрезок\/} 
$[a,b]:=\{x\in \overline\RR\colon a\leq x\leq b\}$ с {\it концами\/}  $a,b\in \overline \RR$, где   $[a,b]=\varnothing$ --- {\it пустое множество\/} при  $a< b$, а также $(a,b]:=[a,b]\setminus a$, $[a,b):=[a,b]\setminus b$ и {\it открытый интервал\/}
$(a,b):=(a,b]\cap [a,b)$.  По определению $\inf \varnothing :=+\infty$ и  $\sup \varnothing :=-\infty$.
{\it Правые и левые открытые полуплоскости\/} обозначаем  как $\CC_{\rh}:= \{z\in \CC \colon \Re z>0\}$ и ${\CC_{\lh}}:=-\CC_{\rh}$.

Для $x\in X\subset \overline \RR$ полагаем $x^+:=\sup\{0,x \}$, где $\sup$ --- точная верхняя грань, $X^+:=\{x^+\colon  x\in X\}$. {\it Расширенной числовой функции\/} $f\colon S\to \overline \RR$ сопоставляем её {\it положительную часть\/} $f^+\colon s\underset{\text{\tiny $s\!\in\! S$}}{\longmapsto} (f(s))^+\in \overline{\RR}^+$ 
 и {\it отрицательную часть\/} $f^-:=(-f)^+\colon S\to \overline{\RR}^+ $.

{\it Распределению точек\/}  ${\sf Z}=\{{\sf z}_{\tt j}\}$ из $\CC$, состоящему из пронумерованных не более чем счётным количеством   индексов ${\tt j}$ точек ${\sf z}_{\tt j}\in \CC$, сопоставляем {\it считающую меру\/}  
 \cite[0.1.2]{Khsur}
\begin{subequations}\label{n}
\begin{align}
n_{\sf Z}&\colon  S\underset{S\subset \CC}{\longmapsto} 
\sum_{{\sf z}_{\tt j}\in S}1\in   \overline \NN_0
\tag{\ref{n}n}\label{df:divmn}\\
\intertext{--- число точек ${\sf z}_{\tt j}$, попавших в $S$. При этом то же обозначение}
n_{\sf Z}&\colon z\underset{\text{\tiny $z\in \CC$}}{\longmapsto} n_{\sf Z}(z)
=\sum_{{\sf z}_{\tt j}=z}1\in  \overline \NN_0
\tag{\ref{n}z}\label{df:divz}\\
\intertext{используется и для {\it считающей функции\/} распределения точек  ${\sf Z}$, а}
n_{\sf Z}^{\rad}(r)&\overset{\eqref{df:divz}}{\underset{r\in \RR^+}{:=}}
\sum_{|z|\leq r}n_{\sf Z}(z)=
\sum_{|{\sf z}_{\tt j}|\leq r}1\in   \overline \NN_0
\tag{\ref{n}r}\label{nrad}
\end{align}
\end{subequations}
--- 
{\it радиальная считающая функция\/} для распределения точек ${\sf Z}$.

Распределения точек    ${\sf Z}$ и  ${\sf Z'}$  {\it совпадают,\/} или {\it равны\/,} и пишем  ${\sf Z}={\sf Z'}$  если у них одни и те же считающие меры или функции $n_{{\sf Z}}\overset{\eqref{df:divmn}}{=}n_{\sf Z'}$, а  включение ${\sf Z}\subset {\sf Z}'$ означает, что  $n_{\sf Z}\leq n_{\sf Z'}$. {\it Объединение\/}  ${\sf Z}\cup {\sf Z'}$ определяется считающей мерой или функцией  $n_{{\sf Z}\cup {\sf Z'}}:=n_{\sf Z}+n_{\sf Z'}$, а при    ${\sf Z}\subset  {\sf Z'}$ {\it разность\/} 
${\sf Z'}\setminus  {\sf Z}$ определяется считающей мерой или функцией  $n_{{\sf Z'}\setminus {\sf Z}}:=n_{\sf Z'}-n_{\sf Z}$. 
Точка $z\in \CC$ {\it принадлежит\/} ${\sf Z}$, т.е. $z\in {\sf Z}$, если $n_{\sf Z}(z)> 0$ для считающей функции \eqref{df:divz}. 
 Распределения точек ${\sf Z}$ {\it конечной верхней плотности\/}  (при порядке $1$), если   конечен верхний предел \cite{Levin56}, \cite{Levin96}
\begin{equation*}
\limsup_{r\to +\infty} \frac{n_{\sf Z}^{\rad}(r)}{r}\in \RR^+. 
\end{equation*}

Ключевую роль в нашей заметке будут играть логарифмические характеристики для распределений точек на $\CC$, введённые для {\it положительных\/} распределений точек ${\sf Z}\subset \RR^+$ в основополагающей для настоящего исследования 
статье   П. Мальявена и Л. А. Рубела \cite{MR} (см. также монографию  Л. А. Рубела \cite{RC}) и распространённые на произвольные {\it комплексные\/} распределения точек   ${\sf Z}\subset \CC$ в работах Б. Н. Хабибуллина 
\cite{KhaD88}, \cite{Kha89}, \cite{kh91AA}, \cite[3.2]{Khsur}. 

Определим {\it  правый и левый характеристические логарифмы} для ${\sf Z}\subset \CC$ как
\begin{subequations}\label{logZC}
\begin{align}
l_{{\sf Z}}^{\rh }(r)&:=\sum_{\substack{0 < |{\sf z}_k|\leq r\\{\sf z}_k \in \CC_{\rh }}} \Re \frac{1}{{\sf z}_k}
=\sum_{0< |{\sf z}_k|\leq r} \Re^+ \frac{1}{{\sf z}_k}, \quad 0<r\leq +\infty, 
\tag{\ref{logZC}r}\label{df:dD+}\\
l_{{\sf Z}}^{\lh }(r)&:=
\sum_{\substack{0< |{\sf z}_k|\leq r\\{\sf z}_k \in \CC_{\lh }}} 
-\Re \frac{1}{{\sf z}_k}=\sum_{0< |{\sf z}_k|\leq r} \Re^- \frac{1}{{\sf z}_k}, \quad 0<r\leq +\infty,
\tag{\ref{logZC}l}\label{df:dD-}
\end{align}
\end{subequations}
а также  {\it правую\/} и {\it левую  логарифмические меры интервалов\/}    $(r,R]\subset \overline{\RR}^+$ как
\begin{subequations}\label{df:l}
\begin{align}
l_{{\sf Z}}^{\rh }(r, R)&\overset{\eqref{df:dD+}}{:=}
l_{{\sf Z}}^{\rh }(R)-l_{{\sf Z}}^{\rh }(r),
\quad 0< r < R \leq +\infty,
\tag{\ref{df:l}r}\label{df:dDl+}\\
l_{{\sf Z}}^{\lh }(r, R)&\overset{\eqref{df:dD-}}{:=}l_{{\sf Z}}^{\lh }(R)-l_{{\sf Z}}^{\lh }(r), \quad 0< r < R \leq +\infty,
\tag{\ref{df:l}l}\label{df:dDl-}
\\
\intertext{которые порождают {\it логарифмическую субмеру  интервалов\/} $(r,R]\subset \overline{\RR}^+$ для  $\sf Z$:} 
l_{{\sf Z}}(r, R)&:=\max \{ l_{{\sf Z}}^{\lh }(r, R), l_{{\sf Z}}^{\rh }(r,
R)\}, \quad 0< r < R \leq +\infty ,
\tag{\ref{df:l}m}\label{df:dDlL}
\end{align}
\end{subequations}
где для ${\sf Z}=\varnothing$ по определению $l_{\varnothing}(r,R)\equiv 0$ при всех $0< r < R \leq +\infty$.

В \cite[(0.2)]{kh91AA} распределение точек ${\sf Z}=\{{\sf z}_{\tt j}\}$ называется {\it отделённым\/}  (углами) {\it от $i\RR$,\/} если   
\begin{equation}\label{dZRe}
|\Re {\sf z}_{\tt j} |\geq d|{\sf z}_{\tt j}| \quad \text{при всех ${\tt j}$ для некоторого числа $d>0$}.
\end{equation} 
Условие \eqref{dZRe} геометрически означает, что  все ненулевые точки из ${\sf Z}$
лежат вне пара непустых открытых вертикальных углов, содержащих $i\RR\setminus 0$. 
 
Распределение точек  ${\sf Z} =\{ {\sf z}_{\tt j}\} \subset \CC$  {\it асимптотически отделено углами  от\/}
$i\RR$,  если 
\begin{equation}\label{con:dis}
\left(\liminf_{{\tt j}\to\infty}\frac{\bigl| \Re {\sf z}_{\tt j} \bigr|}{ |{\sf z}_{\tt j} |} >0 \right)
\underset{\text{или}}{\Longleftrightarrow}
\left(\limsup_{{\tt j}\to\infty}\frac{\bigl| \Im {\sf z}_{\tt j} \bigr|}{ |{\sf z}_{\tt j} |} <1 \right).
\end{equation}
 Пара эквивалентных ограничений \eqref{con:dis} геометрически означает, что  найдётся  пара непустых открытых вертикальных углов, содержащих $i\RR\setminus 0$, для которой точки ${\sf z}_{\tt j}$ лежат вне этой пары углов  при всех ${\tt j}$ за исключением их конечного числа.

Кольцо $\Hol (\CC)$ над $\CC$ состоит из всех голоморфных функций на $\CC$, т.е. $\Hol (\CC)$ --- кольцо {\it целых функций.\/} Через $\Hol_*(\CC):=\bigl\{f\in \Hol(\CC)\colon f\neq 0\bigr\}$ обозначаем множество  всех ненулевых целых функций. Через $\Zero_f$ обозначаем {\it распределение всех  корней\/} целой функции $f\neq 0$ со считающей функцией $n_{\Zero_f}$ в смысле  \eqref{df:divz}, равной  в каждой точке $z\in \CC$ кратности корня функции $f$ в точке $z$.
Целая функция $f\neq 0$ на  $\CC$ \textit{обращается в нуль\/} на распределении  точек ${\sf Z}$ и пишем $f({\sf Z})=0$, если ${\sf Z}\subset {\Zero}_f$.

Целую функцию $f\in \Hol_*(\CC)$ называют \textit{целой функции экспоненциального типа} (пишем {\it ц.ф.э.т}), 
если  для её типа $\type_f$ (при порядке $1$) \cite{Levin56}, \cite{Levin96}, \cite[2.1]{KhI} имеем 
 \begin{equation}\label{typef}
\type_f:=\limsup_{z\to \infty}\frac{\ln |f(z)|}{|z|}<+\infty.
\end{equation} 

Распределение точек   ${\sf Z}=\{{\sf z}_{\tt j}\}\subset \CC$ 
порождает идеал \cite{MR},  \cite[гл. 22]{RC}, \cite{SalKha20}
\begin{equation*}
I({\sf Z}):=\bigl\{f\in \Hol (\CC)\colon f({\sf Z})=0\bigr\}\subset \Hol (\CC)
\end{equation*} 
в кольце $\Hol (\CC)$, а также идеал в кольце всех ц.ф.э.т.
\footnote{В \cite{MR} и  \cite{RC} идеал $I^1({\sf Z})$ обозначен соответственно  как $\mathcal F ({\sf Z})$ и $F({\sf Z})$.}
\begin{equation*}
I^1({\sf Z}):=I({\sf Z})\cap \bigl\{f\in \Hol (\CC)\colon {\type}_f\overset{\eqref{typef}}{<}+\infty \bigr\},
\end{equation*} 
для которых  полагаем 
\begin{equation*}
I_*({\sf Z}):=I({\sf Z})\cap \Hol_*(\CC),
\quad 
I_*^1({\sf Z}):=I^1({\sf Z})\cap \Hol_*(\CC).
\end{equation*}

\begin{thA}[{\cite[теорема 4.1]{MR}, \cite[гл. 22, основная теорема]{RC}}] Пусть  распределения положительных точек  ${\sf Z}\subset \RR^+$ и ${\sf W}\subset \RR^+$   конечной верхней плотности. 
Тогда эквивалентны три утверждения: 
\begin{enumerate}[{\rm I.}]
\item\label{fgi0} Для любой функции $g\in I_*^1({\sf W})$  найдётся функция $f\in I^1_*({\sf Z})$ с ограничением 
\begin{equation}\label{fgiR0}
\bigl|f(iy)\bigr|\leq  \bigl|g(iy)\bigr|\quad \text{при всех  $y\in \RR$}.
\end{equation}
\item\label{fgiii0} Существует пара $f\in I_*^1({\sf Z})$ и  $g\in I^1_*({\sf W})$ c $\Zero_g \cap \CC_{\rh} ={\sf W}$, удовлетворяющая   \eqref{fgiR0}. 
\item\label{fgii0} Существуют $C\in \RR^+$, для которого  
\begin{equation}\label{Zld0}
l_{\sf Z}(r,R)\leq l_{\sf W}(r,R) +C\quad \text{при всех\/ $0< r<R<+\infty$}.
\end{equation}
\end{enumerate}
\end{thA}

В нашей заметке  мы полностью  переносим  теорему Мальявена\,--\,Рубела с положительных на комплексные распределения ${\sf Z}\subset \CC$ и ${\sf W}\subset \CC_{\rh}$, асимптотически отделённые углами от $i\RR$.

\begin{theorem}\label{th1C} Эквивалентность утверждений\/  {\rm  \ref{fgi0}--\ref{fgii0}} теоремы Мальявена\,--\,Рубела имеет место и для  любых распределений комплексных точек ${\sf Z}\subset \CC$ и ${\sf W}\subset \CC_{\rh}$ конечной верхней плотности, асимптотически отделённых углами от мнимой оси $i\RR$ в смысле \eqref{con:dis}.
\end{theorem}
\begin{proof} Для доказательства импликации \ref{fgi0}$\Longrightarrow$\ref{fgiii0} достаточно построить ц.ф.э.т. $g\neq 0$ с $\Zero_g \cap \CC_{\rh} ={\sf W}$. В качестве такой функции можно выбрать ц.ф.э.т. $g$ в виде чётного канонического произведения Адамара\,--\,Вейерштрасса    \cite{Levin56}, \cite{Levin96}
\begin{equation*}
g(z)\underset{z\in \CC}{:=}\prod_{{\tt j}}\Bigl(1-\frac{z^2}{{\sf w}_{\tt j}^2}\Bigr). 
\end{equation*}
Для доказательства импликации \ref{fgiii0}$\Longrightarrow$\ref{fgii0}  в обозначении 
\begin{equation*}
J_{i\RR}(r,R;v):=
\frac{1}{2\pi}\int_r^{R} \frac{v(-iy)+v(iy)}{y^2} \dd y, \quad 0<r<R<+\infty.
\end{equation*} 
для борелевских  функций $v\colon i\RR \to \overline \RR$ будет  использована 
\begin{lemma}[{\cite[(1.3)]{Kha89}, \cite[(0.4)]{kh91AA} и в явном  в \cite[предложение 4.1, (4.19)]{KhII}}]\label{lemJl} При любом фиксированном числе $r_0>0$ для любой ц.ф.э.т. $f\neq 0$ имеет место соотношение
\begin{equation}\label{Jll}
\sup_{r_0\leq r<R<+\infty}\max \Bigl\{\bigl|J_{i\RR}(r,R;\ln |f|)-l_{\Zero_f}^{\rh}(r,R)\bigr|,\;
\bigl|J_{i\RR}(r,R;\ln|f|)-l_{\Zero_f}^{\lh}(r,R)\bigr|\Bigr\}
<+\infty.
\end{equation}
\end{lemma}
В условиях утверждения \ref{fgiii0} логарифмирование и интегрирование  неравенства \eqref{fgiR0} даёт 
\begin{equation}\label{Jfg}
J_{i\RR}\bigl(r,R;\ln |f|\bigr)\leq J_{i\RR}\bigl(r,R;\ln |g|\bigr)\quad\text{для всех $r_0\leq r<R<+\infty$}.
\end{equation} 
По лемме \ref{lemJl}, применённой  к ц.ф.э.т. $f\neq 0$, 
для некоторого числа $C_f\in \RR^+$ получаем
\begin{equation*}\label{Zdg}
l_{\sf Z}(r,R)\overset{\eqref{df:dDlL}}{\leq} l_{\Zero_f}(r,R)\overset{\eqref{Jll}}{\leq} 
J_{i\RR}\bigl(r,R;\ln |f|\bigr) +C_f\overset{\eqref{Jfg}}{\leq}
 J_{i\RR}\bigl(r,R;\ln |g|\bigr)+C_f
\end{equation*}
для всех $r_0\leq r<R<+\infty$. Вновь по  лемме \ref{lemJl}, но применённой  уже  к ц.ф.э.т. $g\neq 0$,  можем продолжить эту цепочку неравенств с некоторым числом $C_g\in \RR^+$ как  
$$
 l_{\sf Z}(r,R)\leq  J_{i\RR}\bigl(r,R;\ln |g|\bigr)+C_f\overset{\eqref{Jll}}{\leq} 
l_{\Zero_f}^{\rh}(r,R)+C_g+C_f=
l_{\sf W}^{\rh}(r,R) +C_g+C_f=l_{\sf W}(r,R) +C_g+C_f
$$
при всех\/ $r_0\leq  r<R<+\infty$, где последние равенства следуют из  определения 
\eqref{df:dD+} правой логарифмической меры и из условия $\Zero_g \cap \CC_{\rh} ={\sf W}\subset \CC_{\lh }$ утверждения \ref{fgiii0}. Отсюда для постоянной $C:=C_g+C_f$, не зависящей от $r>r_0$ и $R>r$, по определениям \eqref{df:l}
при достаточно малом  $r_0>0$ можно рассматривать все  $0<r<R<+\infty$, что даёт требуемые неравенства  \eqref{Zld0} и соответственно утверждение \ref{fgii0}.

Для доказательства  импликации \ref{fgii0}$\Longrightarrow$\ref{fgi0} 
будет использована 
\begin{thD}[{\cite[основная теорема]{kh91AA}}] Пусть $g\neq 0$ --- ц.ф.э.т., а  распределение корней $\Zero_g$ и распределение точек ${\sf Z}=\{{\sf z}_{\tt j}\}$ отделены (углами) от $i\RR$ в смысле \eqref{dZRe}. 
 Для существования ц.ф.э.т. $f\neq 0$, обращающейся в нуль на ${\sf Z}$ и удовлетворяющей  \eqref{fgiR0}, необходимо и достаточно, чтобы существовало число $M\in \RR$ такое, что 
\begin{equation}\label{fgM}
l_{\sf Z}(r,R) \leq J_{i\RR}\bigl(r,R;\ln |g|\bigr)+M\quad\text{для всех $1\leq r<R<+\infty.$}
\end{equation}
 \end{thD}
Пусть в условиях утверждения \ref{fgii0}    выполнены условия теоремы Kh. Тогда для ц.ф.э.т. $g\neq 0$ с ${\sf W}\subset \Zero_g$ по лемме \ref{lemJl} существует число $C_g\in \RR^+$, с которым 
$$
l_{\sf Z}(r,R)\overset{\eqref{Zld0}}{\leq} l_{\sf W}(r,R) +C
\leq l_{\Zero_g}(r,R) +C\overset{\eqref{Jll}}{\leq}  
 J_{i\RR}\bigl(r,R;\ln |g|\bigr)+C_g+C 
$$
при всех $r_0\leq r<R<+\infty$. Отсюда при $r_0:=1$ получаем \eqref{fgM} и из части достаточности в теореме Kh найдётся
ц.ф.э.т. $f\neq 0$ с $f({\sf Z})=0$, удовлетворяющая   \eqref{fgiR0}. 

Пусть теперь распределения точек ${\sf Z}$ и  $\Zero_g\supset {\sf W}$ лишь асимптотически отделены углами от $i\RR$ в смысле \eqref{con:dis}.  Всегда можно выделить конечное распределение точек ${\sf Z}_0\subset {\sf Z}$ и конечное распределение точек ${\sf G}_0\subset \Zero_g$ так, что  ${\sf Z}_{\infty}:= {\sf Z}\setminus {\sf Z}_0$ и ${\sf G}_{\infty}:= {\Zero_g}\setminus {\sf G}_0$ уже отделены углами от $i\RR$ в смысле \eqref{dZRe}. При этом ввиду конечности распределения точек  ${\sf G}_0$ по определениям логарифмических функций интервалов \eqref{logZC}--\eqref{df:l} найдётся число $C_0\in \RR^+$, для которого $l_{{\sf G}_0}(r,R)\leq C_0$ при всех $0<r<R<+\infty$, откуда 
\begin{multline}\label{lW}
l_{{\sf Z}_{\infty}}(r,R)\leq l_{{\sf Z}}(r,R)\overset{\eqref{Zld0}}{\leq} l_{\sf W}(r,R) +C\leq 
l_{\Zero_g}(r,R) +C=l_{{\sf G}_{\infty}\cup {\sf G}_{0}}(r,R) +C
\\
\overset{\eqref{df:dDlL}}{\leq} l_{{\sf G}_{\infty}}(r,R) + l_{{\sf G}_{0}}(r,R)+C
\leq l_{{\sf G}_{\infty}}(r,R) +C_0+C \quad\text{для всех $0<r<R<+\infty$}.
\end{multline}
Рассмотрим теперь ц.ф.э.т. $g_{\infty}:=g/g_0\neq 0$, где $g_0$ --- некоторый многочлен с распределением  корней  
${\sf G}_0$, с  $g_{\infty}({\sf G}_{\infty})=0$. По уже доказанной  версии импликации \ref{fgii0}$\Longrightarrow$\ref{fgi0} 
для  $g_{\infty}$ в роли $g$ и ${\sf G}_{\infty}$ в роли ${\sf W}$ ввиду  \eqref{lW}  найдётся ц.ф.э.т. $f_{\infty}\neq 0$ 
с $f_{\infty}({\sf Z}_{\infty})=0$, удовлетворяющая условию \eqref{fgiR0} в виде 
$\bigl|f_{\infty}(iy)\bigr|\leq \bigl|g_{\infty}(iy)\bigr|$ при всех  $y\in \RR$, откуда  по построениям
\begin{equation}\label{gfW}
 \bigl|(g_0f_{\infty})(iy)\bigr|\leq \bigl|g_0g_{\infty}(iy)\bigr|=\bigl|g(iy)\bigr|
 \quad\text{при всех  $y\in \RR$ и} \quad (g_0f_{\infty})({\sf Z}_{\infty})=0.
\end{equation}
Пусть число точек в ${\sf Z}_0$ равно $N$, т.е. $N:=n_{{\sf Z}_0}(\CC)$,  и $f_0$ --- некоторый многочлен с распределением степени $N$ корней ${\sf Z_0}$. Тогда для достаточно малого числа $a>0$ сужение функции 
\begin{equation*}
f_a(z)\underset{z\in \CC}{:=}af_0(z)\Bigl(\frac{\sin iz}{z}\Bigr)^N 
\end{equation*}
на $i\RR$ ограничено по модулю единицей на $i\RR$  и по  \eqref{gfW} для ц.ф.э.т. $f:=f_ag_0f_{\infty}\neq 0$, обращающейся в нуль  на ${\sf Z}={\sf Z}_0\cup {\sf Z}_{\infty}$ имеем  
$ |f(iy)|= \bigl|(f_ag_0f_{\infty})(iy)\bigr|\overset{\eqref{gfW}}{\leq} \bigl|g(iy)\bigr|$
при всех  $y\in \RR$, что завершает доказательство импликации \ref{fgii0}$\Longrightarrow$\ref{fgi0} и теоремы \ref{th1C}. 
\end{proof}
\begin{remark}\label{rem1}
Анализ доказательства теоремы \ref{th1C} показывает, что при доказательстве импликации   
\ref{fgii0}$\Longrightarrow$\ref{fgi0} нигде не используется расположение распределения точек ${\sf W}$ именно в правой полуплоскости $\CC_{\rh}$, а  условие    ${\sf W}\subset \CC_{\rh}$ в теореме \ref{th1C} применяется только при доказательстве импликаций \ref{fgi0}$\Longrightarrow$\ref{fgiii0}$\Longrightarrow$\ref{fgii0}. Ясно, что условие  
${\sf W}\subset \CC_{\rh}$ можно заменить на расположение ${\sf W}\subset \CC_{\lh }$ в левой полуплоскости с помощью зеркальной симметрии относительно мнимой оси. 
\end{remark}

\end{document}